\documentclass{elsart}
\usepackage{amsmath,amssymb,epsfig}

 \newcounter{mmacnt}
 \def\restartmma{\setcounter{mmacnt}{0}}
 \restartmma \catcode`|=\active
 \def|#1|{\mathrm{#1}}
 \catcode`|=12

\begin{document}

\begin{frontmatter}

  \title{Inequalities for ultraspherical polynomials.\\
   Proof of a conjecture of I. Ra\c{s}a}

\author{Geno Nikolov}
\address{Department of Mathematics and Informatics,
         Sofia University ``St. Kliment Ohridski'' \\
         5 James Bourchier Blvd., 1164 Sofia, Bulgaria}
\ead{geno@fmi.uni-sofia.bg}

\begin{abstract}
A recent conjecture by I. Ra\c{s}a asserts that the sum of the
squared Bernstein basis polynomials is a convex function in $[0,1]$.
This conjecture turns out to be equivalent to a certain upper
pointwise estimate of the ratio $P_n^{\prime}(x)/P_n(x)$ for $x\geq
1$, where $P_n$ is the $n$-th Legendre polynomial. Here, we prove
both upper and lower pointwise estimates for the ratios
$\big(P_n^{(\lambda)}(x)\big)^{\prime}/P_n^{(\lambda)}(x)$, $~x\geq
1$, where $P_n^{(\lambda)}$ is the $n$-th ultraspherical polynomial.
In particular, we validate Ra\c{s}a's conjecture.
\end{abstract}

\end{frontmatter}

\section{Introduction and statement of the results}

The classical Bernstein operator
\[
B_n(f;t)=\sum_{k=0}^{n}b_{n,k}(t)f\Big(\frac{k}{n}\Big),\ \
\]
where the basis polynomials $b_{n,k}$ are given by
\[
b_{n,k}(t)=\binom{n}{k}t^k(1-t)^{n-k},\ \ k=0,1,\ldots,n,
\]
plays an important role in Approximation Theory, and has been an
object of intensive study throughout the years, see  \cite[Chapter
10]{DL:1993}.

Recently, Ioan Ra\c{s}a \cite{IR:2012} (see \cite{GRR:2014} for a
more accessible source) formulated the following conjecture about
the sum of the squared Bernstein basis polynomials:

\begin{conj}\label{con}
The function
\[
F_n(t):= \sum_{k=0}^{n}b_{n,k}^2(t)
\]
is convex in $[0,1]$.
\end{conj}

As a matter of fact, Conjecture \ref{con} is the strongest amongst
three conjectures proposed by Ra\c{s}a in \cite{IR:2012}; the other
two are that $F_n(t)$ attains its minimum in $[0,1]$ at $t=1/2$, and
that $F_n(t))$ is monotone decreasing in $[0,1/2]$ and monotone
increasing in $[1/2,1]$ (see \cite[Conjectures 3.4 and
3.3]{GRR:2014}). These weaker conjectures have been validated by
Thorsten Neuschel \cite{TN:2013} (see also
\cite[Lemma~3.5]{GRR:2014}). The function $F_n$ appears as a factor
in an upper estimate of the degree of non-multiplicativity of
Bernstein operators, see \cite{GRR:2014}.

In Section 2 we prove the following equivalence:

\begin{thm}\label{thm:1}
Conjecture \ref{con} is equivalent to the inequality
\begin{equation}\label{e1.1}
\frac{P_n^{\prime}(x)}{P_n(x)}\leq\frac{2n^2}{x+(2n-1)\sqrt{x^2-1}},\
\ x\geq 1\,,
\end{equation}
where $P_n$ is the $n$-th Legendre polynomial.
\end{thm}
As is well-known, the Legendre polynomial $P_n$ belongs to the
family of ultraspherical polynomials $P_n^{(\lambda)}$,
$\lambda>-1/2$, which are orthogonal in $[-1,1]$ with respect to the
weight function $w_{\lambda}(x)=(1-x^2)^{\lambda-1/2}$ ($P_n$
corresponds to the case $\lambda=1/2$). This is our motivation
instead of proving inequality \eqref{e1.1} to obtain upper (but also
lower) pointwise bounds for the ratios
$\big(P_n^{(\lambda)}(x)\big)^{\prime}/P_n^{(\lambda)}(x)$, $~x\geq
1$. For the sake of brevity, hereafter we skip the superscript
$^{(\lambda)}$ and write
\[
p_n(x):=P_n^{(\lambda)}(x)\,,\ \ \ \
u_n(x):=\frac{p_n^{\prime}(x)}{p_n(x)}\,.
\]
In Section 3 we prove

\begin{thm}\label{thm:2} Let $n\in \mathbb{N}$ and $\lambda
>-1/2$. Then
\begin{equation}\label{e1.2}
u_n(x)\geq\frac{n(n+2\lambda)} {(2\lambda+1)x+(n-1)\sqrt{x^2-1}}\,,\
\ x\geq 1\,.
\end{equation}
Moreover, if $\lambda\in [0,1]$, then
\begin{equation}\label{e1.3}
u_n(x)\leq\frac{n^2(n+\lambda)}
{\lambda(n+1)x+(n^2-\lambda)\sqrt{x^2-1}}\,,\ \ x\geq 1\,.
\end{equation}
\end{thm}

Specialized to the case $\lambda=1/2$, inequality \eqref{e1.3}
becomes
\[
\frac{P_n^{\prime}(x)}{P_n(x)}\leq \frac{n^2(2n+1)}
{(n+1)x+(2n^2-1)\sqrt{x^2-1}}\,,\ \ x\geq 1\,,
\]
which is easily seen to be slightly stronger than \eqref{e1.1}.
Consequently, we have
\begin{cor}\label{cor:1}
Conjecture \ref{con} is true.
\end{cor}

Yet another theorem of the same nature is

\begin{thm}\label{thm:3} Let $n\in \mathbb{N}$ and $\lambda
>-1/2$. Then
\begin{equation}\label{e1.4}
n\Big(\frac{1}{x}+\frac{n-1}{2(n+\lambda-1)x^3}\Big)\leq u_n(x) \leq
n\Big(\frac{1}{x}+\frac{n-1}{(2\lambda+1)x^3}\Big)\,,\ \ \ x\geq 1.
\end{equation}
\end{thm}
The proof of Theorem \ref{thm:3} is given in Section 4. In the final
section we compare Theorems \ref{thm:2} and \ref{thm:3} and discuss
their sharpness.

\section{Equivalent formulation of Conjecture 1: Proof of Theorem \ref{thm:1}}
\setcounter{equation}{0} Clearly,
\[
F_n(t)=\sum_{k=0}^{n}\binom{n}{k}^2t^{2k}(1-t)^{2n-2k}
\]
is a polynomial of degree $2n$, and $F_n(t)=F_n(1-t)$. Therefore,
Conjecture \ref{con} is equivalent to the inequality
\begin{equation}\label{e2.1}
F_n^{\prime\prime}(t)\geq 0,\ \  t\in [0,1/2).
\end{equation}
We shall express $F_n$ through the Legendre polynomial $P_n$. By
Rodrigues' formula and the Leibnitz rule we have
\begin{equation}\label{e2.2}
P_n(x)=\frac{1}{2^n\,n!}\,\frac{d^{n}}{dx^{n}}\big\{(x^2-1)^{n}\big\}
=\frac{1}{2^{n}}\sum_{k=0}^{n}\binom{n}{k}^2(x-1)^k(x+1)^{n-k}\,.
\end{equation}
Let us set
\[
x=x(t):=\frac{1}{2}\Big(1-2t+\frac{1}{1-2t}\Big),\ \ t\in [0,1/2),
\]
then $x$ traces the interval $[1,\infty)$, and
\[
x-1=\frac{2t^2}{1-2t},\ \ x+1=\frac{2(1-t)^2}{1-2t},\
1-2t=x-\sqrt{x^2-1}\,.
\]
Replacement in \eqref{e2.2} yields
\[
P_n(x)=\frac{1}{(1-2t)^{n}}\sum_{k=0}^{n}\binom{n}{k}^2t^{2k}(1-t)^{2n-2k}
=(1-2t)^{-n}F_n(t)\,,
\]
hence
\begin{equation}\label{e2.3}
F_n(t)=\big(x-\sqrt{x^2-1}\big)^n P_n(x)\,.
\end{equation}
Taking into account that
\[
x^{\prime}(t)=\frac{1}{(1-2t)^2}-1=\frac{4t(1-t)}{(1-2t)^2}=
\frac{2\sqrt{x^2-1}}{x-\sqrt{x^2-1}},
\]
we differentiate \eqref{e2.3} to obtain consecutively
\begin{equation}\label{e2.4}
\begin{split}
F_n^{\prime}(t)&=\frac{d}{dx}\,\big\{(x-\sqrt{x^2-1})^n
P_n(x)\big\}\,x^{\prime}(t)\\
&=2(x-\sqrt{x^2-1})^{n-1} \big(\sqrt{x^2-1}\,P_n^{\prime}(x)-n
P_n(x)\big)
\end{split}
\end{equation}
(this formula has been  already obtained by Neuschel
\cite{TN:2013}),
\[
\begin{split}
F_n^{\prime\prime}(t)&=2\frac{d}{dx}\,\big\{(x-\sqrt{x^2-1})^{n-1}
\big(\sqrt{x^2-1}\,P_n^{\prime}(x)-n
P_n(x)\big)\big\}\,x^{\prime}(t)\\
&=4(x-\sqrt{x^2-1})^{n-2}\\
&\quad\times\Big((x^2-1)\,P_n^{\prime\prime}(x)+\big(x-(2n-1)\sqrt{x^2-1}\big)P_n^{\prime}(x)
+n(n-1)P_n(x)\Big)\,.
\end{split}
\]
On using the differential equation
\[
(1-x^2)P_n^{\prime\prime}(x)-2x\,P_n^{\prime}(x)+n(n+1)P_n(x)=0\,,
\]
we replace $(x^2-1)P_n^{\prime\prime}(x)$ to finally obtain
\begin{equation}\label{e2.5}
F_n^{\prime\prime}(t)=4(x-\sqrt{x^2-1})^{n-2}\Big[2n^2P_n(x)-
\big(x+(2n-1)\sqrt{x^2-1}\big)P_n^{\prime}(x)\Big]\,.
\end{equation}
Since $x-\sqrt{x^2-1}>0$ for $x\geq 1$, it follows that \eqref{e2.1}
is fulfilled exactly when the term in the square brackets in
\eqref{e2.5} is non-negative. The latter is equivalent to
\eqref{e1.1}, as $P_n(x)>0$ for $x\geq 1$. Theorem \ref{thm:1} is
proved.
\section{Proof of Theorem \ref{thm:2}}
\setcounter{equation}{0} For easy reference, we collect in a lemma
some well-known properties of ultraspherical polynomials, see, e.g.,
\cite[Chapter 4.7]{GS:1959}.

\begin{lem}\label{lem:1}
For $\lambda\ne 0$, $p_n=P_n^{(\lambda)}$ satisfies the following
properties:

{\rm ~~(i)~~~~}
$\displaystyle{(1-x^2)p_n^{\prime\prime}(x)-(2\lambda+1)x\,p_n^{\prime}(x)
+n(n+2\lambda)p_n(x)=0}$;

{\rm ~~(ii)~~~~} $\displaystyle{p_{n+1}^{\prime}(x)=
(n+2\lambda)p_n(x)+x\,p_{n}^{\prime}(x)}$;

{\rm ~~(iii)~~~~} $\displaystyle{p_{n+1}(x)=
\frac{1}{n+1}\,\big((n+2\lambda)x\,p_n(x)+(x^2-1)\,p_{n}^{\prime}(x)\big)}$;

{\rm ~~(iv)~~~~} $\displaystyle{n\,p_{n}(x)=
x\,p_n^{\prime}(x)-p_{n-1}^{\prime}(x)}$;

{\rm ~~(v)~~~~} $\displaystyle{(n+1)p_{n+1}(x)=
2(n+\lambda)x\,p_n(x)-(n+2\lambda-1)\,p_{n-1}(x)}$,~~~$n\ge 1$.
\end{lem}

Set
\[
u_n(x):=\frac{p_n^{\prime}(x)}{p_n(x)},\ \ n\in \mathbb{N},
\]
then obviously $u_n(x)$ is positive and strictly monotone decreasing
in $[1,\infty)$. On using Lemma \ref{lem:1} (i), we find
\begin{equation}\label{e3.1}
u_n(1)=\frac{n(n+2\lambda)}{2\lambda+1}\,.
\end{equation}

The proof of Theorem \ref{thm:2} goes by induction. For the
induction transition from $n$ to $n+1$, we make use of Lemma
\ref{lem:1} (ii) and (iii) to obtain
\begin{equation}\label{e3.2}
u_{n+1}(x)=(n+1)\, \frac{n+2\lambda+x\,u_n(x)}
{(n+2\lambda)x+(x^2-1)\,u_n(x)}\,.
\end{equation}
We observe that the function
$\varphi(t)=\frac{a+x\,t}{a\,x+(x^2-1)t}$ is continuous and strictly
monotone increasing in $(0,\infty)$ whenever $a>0$ and $x\geq 1$.

Let us prove first inequality \eqref{e1.2}. Clearly $u_1(x)=1/x$
satisfies \eqref{e1.2} with equality for every $x\geq 1$. Suppose
that, for some $n\in \mathbb{N}$,
\[
u_n(x)\geq \frac{n(n+2\lambda)}
{(2\lambda+1)x+(n-1)\sqrt{x^2-1}}=:t_n(x),\ \ x\geq 1.
\]
Then, by \eqref{e3.2},
\[
u_{n+1}(x)\geq (n+1)\, \frac{n+2\lambda+x\,t_n(x)}
{(n+2\lambda)x+(x^2-1)\,t_n(x)}\,,
\]
and the induction step will be performed once we show that
\[
(n+1)\, \frac{n+2\lambda+x\,t_n(x)}
{(n+2\lambda)x+(x^2-1)\,t_n(x)}\geq t_{n+1}(x),\ \ x\geq 1\,.
\]
A straightforward calculation shows that the latter inequality is
equivalent to the obvious inequality
\[
2(\lambda+1)n\sqrt{x^2-1}\big(x-\sqrt{x^2-1}\big)\geq 0,\ \ x\geq 1,
\]
and this accomplishes the induction proof of \eqref{e1.2}.

Inequality \eqref{e1.3} can be proved by induction along the same
lines as \eqref{e1.2}. However, we would like to provide some clue
about the way we deduced this inequality.

We seek for which non-negative $c_n=c_n(\lambda)$ the inequality
\begin{equation}\label{e3.3}
u_n(x)\leq\frac{n^2}{c_n x+(n-c_n)\sqrt{x^2-1}}=:\tau(n,c_n,x),\ \ \
x\geq 1,
\end{equation}
holds true. As is easy to see, the larger $c_n$, the better (i.e.,
smaller) the upper bound $\tau(n,c_n,x)$ in \eqref{e3.3}. However,
$c_n$ cannot be arbitrarily large, for, according to \eqref{e3.1},
in order that \eqref{e3.3} is true for $x=1$, there must hold
$n^2/c_n\geq n(n+2\lambda)/(2\lambda+1)$. Hence, a natural
restriction for $c_n$ is
\[
0\leq c_n\leq \frac{(2\lambda+1)n}{n+2\lambda}.
\]
Assume that \eqref{e3.3} holds true for some $n\in \mathbb{N}$. By
\eqref{e3.2}, we have
\[
u_{n+1}(x)\leq (n+1)\, \frac{n+2\lambda+x\,\tau(n,c_n,x)}
{(n+2\lambda)x+(x^2-1)\,\tau(n,c_n,x)}\,,\ \ x\geq 1.
\]
The induction step will be done if we manage to show that
\begin{equation}\label{e3.4}
(n+1)\, \frac{n+2\lambda+x\,\tau(n,c_n,x)}
{(n+2\lambda)x+(x^2-1)\,\tau(n,c_n,x)}\leq \tau(n+1,c_{n+1},x),\ \
x\geq 1.
\end{equation}
At this point we assume that the sequence $\{c_n\}_{n=1}^{\infty}$
is non-increasing. Then \eqref{e3.4} will be a consequence of the
inequality
\begin{equation}\label{e3.5}
(n+1)\, \frac{n+2\lambda+x\,\tau(n,c_n,x)}
{(n+2\lambda)x+(x^2-1)\,\tau(n,c_n,x)}\leq \tau(n+1,c_{n},x),\ \
x\geq 1,
\end{equation}
since $\tau(n,c,x)$ is a decreasing function of $c$. Now we check
for which $c_n$ the inequality \eqref{e3.5} holds true. Inequality
\eqref{e3.5} is equivalent (for brevity, here we write $c$ instead
of $c_n$) to
\[
\frac{\big[n^2\!+\!c(n\!+\!2\lambda)\big]x\!+\!(n\!-\!c)
(n\!+\!2\lambda)\sqrt{x^2\!-\!1}}
{c(n\!+\!2\lambda)x^2\!+\!(n\!-\!c)(n\!+\!2\lambda)x\sqrt{x^2\!-\!1}
\!+\!n^2(x^2\!-\!1)}\!\leq\!
\frac{n\!+\!1}{c\,x\!+\!(n\!+\!1\!-\!c)\sqrt{x^2\!-\!1}}.
\]
Since $x\geq 1$, both denominators are positive, and after
simplification and cancelation of the positive factor
$x-\sqrt{x^2-1}$, the above inequality is reduced to the inequality
\[
c\big[2\lambda(n\!+\!1)\!+\!n \!-\!c(n+2\lambda)\big]x
\!-\!\big[(2n\!+\!1)(n\!+\!2\lambda)c\!-\!(n\!+\!2\lambda)c^2\!-\!2\lambda
n(n\!+\!1)\big]\sqrt{x^2\!-\!1}\geq 0.
\]
For the last inequality to be true for every $x\geq 1$, the
coefficient of $x$ must be non-negative and greater than or equal to
the coefficient of $\sqrt{x^2-1}$, i.e., there must hold
\[
c\big[2\lambda(n\!+\!1)\!+\!n \!-\!c(n+2\lambda)\big]\geq \max\{0,
(2n\!+\!1)(n\!+\!2\lambda)c\!-\!(n\!+\!2\lambda)c^2\!-\!2\lambda
n(n\!+\!1)\}.
\]
The latter is equivalent to
\begin{equation}\label{e3.6}
0\leq
c\leq\min\Big\{\frac{2\lambda(n+1)+n}{n+2\lambda},\frac{\lambda(n+1)}{n+\lambda}
\Big\}.
\end{equation}
Clearly, if $-1/2<\lambda<0$, then \eqref{e3.6} to has no solution.
On the other hand, if $0\leq\lambda\leq 1$, then
\[
c=c_n=\frac{\lambda(n+1)}{n+\lambda}
\]
is a solution of \eqref{e3.6}, and the sequence
$\{c_n\}_{n=1}^{\infty}$ is non-increasing, in accordance with our
assumption.

Performing our reasoning backward, we see that if
$c_n=\frac{\lambda(n+1)}{n+\lambda}$, $0\leq \lambda\leq 1$, then
$u_n(x)\leq \tau(n,c_n,x)$ for every $x\geq 1$ implies
$u_{n+1}(x)\leq \tau(n+1,c_{n+1},x)$ for every $x\geq 1$, i.e., the
induction step is done. Notice that for this choice of $c_n$ we have
\[
\tau(n,c_n,x)=\frac{n^2(n+\lambda)}
{\lambda(n+1)x+(n^2-\lambda)\sqrt{x^2-1}},
\]
therefore \eqref{e3.3} is in fact inequality \eqref{e1.3}. It
remains to verify \eqref{e1.3} for $n=1$, i.e.,
\[
\frac{1}{x}\leq \frac{1+\lambda}{2\lambda
x+(1-\lambda)\sqrt{x^2-1}},\ \ \ x\geq 1.
\]
The latter inequality is equivalent to
$(1-\lambda)\big(x-\sqrt{x^2-1}\big)\geq 0$, hence is true. The
proof of Theorem \ref{thm:2} is complete.
\section{Proof of Theorem \ref{thm:3}}
\setcounter{equation}{0} On using Lemma \ref{lem:1} (iv) we obtain
\begin{equation}\label{e4.1}
u_n(x)=\frac{1}{x}\Big(n+\frac{p_{n-1}^{\prime}(x)}{p_n(x)} \Big)\,.
\end{equation}
Let $\{x_k\}_{k=1}^{n}$ be the zeros of $p_n=P_n^{(\lambda)}$, they
form a symmetrical set with respect to the origin. Invoking again
Lemma \ref{lem:1} (iv) we get
\[
\frac{p_{n-1}^{\prime}(x_k)}{p_n^{\prime}(x_k)}=x_k,\ \
k=1,\ldots,n.
\]
By Lagrange's interpolation formula and the symmetry of the set
$\{x_k\}_{k=1}^{n}$ we obtain
\begin{equation}\label{e4.2}
\begin{split}
\frac{p_{n-1}^{\prime}(x)}{p_n(x)}&=\sum_{k=1}^{n}
\frac{p_{n-1}^{\prime}(x_k)}{p_n^{\prime}(x_k)}\cdot\frac{1}{x-x_k}
=\sum_{k=1}^{n}\frac{x_k}{x-x_k}\\
&=\frac{1}{2}\sum_{k=1}^{n}\Big(\frac{x_k}{x-x_k}-\frac{x_k}{x+x_k}\Big)
=\sum_{k=1}^{n}\frac{x_k^2}{x^2-x_k^2}\,.
\end{split}
\end{equation}
Hence,
\[
\psi(x):=\frac{x^2
p_{n-1}^{\prime}(x)}{p_n(x)}=\sum_{k=1}^{n}\frac{x_k^2}{1-(x_k/x)^2}\,,
\]
and
\[
\psi^{\prime}(x)=-\frac{2}{x^3}
\sum_{k=1}^{n}\frac{x_k^4}{\big(1-(x_k/x)^2\big)^2}<0\ \ \mbox{ for
}\ x\geq 1\,.
\]
We observe that $\psi$ is a monotone decreasing function in
$[1,\infty)$, therefore
$\psi(1)\geq\psi(x)\geq\lim_{x\rightarrow\infty}\psi(x)$ therein.
Lemma \ref{lem:1} (iv) and \eqref{e3.1} imply
\begin{equation}\label{e4.3}
\psi(1)=u_n(1)-n=\frac{n(n-1)}{2\lambda+1}\,.
\end{equation}
On the other hand, we have
\[
\lim_{x\rightarrow\infty}\psi(x)=(n-1)\,\frac{a_{n-1}}{a_n}
\]
with $a_m=a_m(\lambda)$ being the leading coefficient of $p_m$. From
$a_0=1$, $a_1=2\lambda$ and Lemma \ref{lem:1} (v) we infer
\[
a_m=\frac{2^m\,\lambda(\lambda+1)\cdots(\lambda+m-1)}{m!},\ \ m\in
\mathbb{N},
\]
whence
\begin{equation}\label{e4.4}
\lim_{x\rightarrow\infty}\psi(x)=\frac{n(n-1)}{2(n+\lambda-1)}\,.
\end{equation}
Thus,
\[
\frac{n(n-1)}{2(n+\lambda-1)x^2}\leq\frac{p_{n-1}^{\prime}(x)}{p_n(x)}
\leq\frac{n(n-1)}{(2\lambda+1)x^2}\,,\ \ x\in [1,\infty)\,.
\]
Theorem \ref{thm:3} follows from substituting these bounds in
\eqref{e4.1}.
\section{Final remarks}
\setcounter{equation}{0} The bounds for
$u_n(x)=p_n^{\prime}(x)/p_n(x)$ provided by Theorems \ref{thm:2} and
\ref{thm:3} are sharp as $x\rightarrow\infty$ in the sense that they
preserve the property $\lim_{x\rightarrow\infty}x\,u_n(x)=n$. The
lower bound in Theorem \ref{thm:2} and the upper bound in Theorem
\ref{thm:3} are also sharp for $x=1$, see \eqref{e3.1}. However,
except for some neighborhoods of $x=1$, the latter bounds are
inferior to their counterparts given in Theorem \ref{thm:3} and
Theorem \ref{thm:2}, respectively.

Unfortunately, the upper bound in Theorem \ref{thm:2} was only
proven for $0\leq\lambda\leq 1$. This restriction on $\lambda$ in
not a proof defect, as for negative $\lambda$ the right-hand side of
\eqref{e1.3} is negative at $x=1$ while if $\lambda>1$ then
\eqref{e1.3} fails for $n=1$.

As is well-known (see, e.g., \cite[Theorem 6.21.1]{GS:1959}), the
squared $k$-th zero $x_k^2=x_k^2(\lambda)$ of $p_n=P_n^{(\lambda)}$
is a monotone decreasing function of $\lambda$, ($k=1,\ldots,n$). By
\eqref{e4.2} we deduce that
\[
\frac{p_{n-1}^{\prime}(x)}{p_n(x)}=\sum_{k=1}^{n}\frac{1}{x^2/x_k^2-1}
\]
is also a monotone decreasing function of $\lambda$ for $x\geq 1$,
and so is $u_n(x)$, by virtue of \eqref{e4.1}. Therefore, the upper
bound for $u_n(x)$ given by \eqref{e1.3} for $\lambda=1$ is also
upper bound for $u_n(x)$ whenever $\lambda\geq 1$, i.e.,
\begin{equation}\label{e5.1}
u_n(x)\leq \frac{n^2} { x+(n-1)\sqrt{x^2-1}} \ \ \mbox{ for every }\
x\geq 1\ \ \mbox{ and }\ \lambda\geq 1\,.
\end{equation}

For particular $\lambda\geq 1$, $n$ and $x>1$ \eqref{e5.1} may
provide better upper bound for $u_n(x)$ than the one given by
Theorem \ref{thm:3}.

Of course, in the cases $\lambda=0$ and $\lambda=1$ one can obtain
explicit formulae for $u_n(x)$ exploiting the representation of the
Chebyshev polynomials of the first and second kind (see, e.g.,
\cite[Eqns. (1.49) and (1.52)]{MH:2000})
\[
T_n(x)=\frac{1}{2}\big((x+\sqrt{x^2-1})^n+(x-\sqrt{x^2-1})^n\big)\,,
\]
\[
U_n(x)=\frac{1}{2\sqrt{x^2-1}}\big((x+\sqrt{x^2-1})^{n+1}-
(x-\sqrt{x^2-1})^{n+1}\big)\,.
\]

Let us finally point out that for fixed $n\in \mathbb{N}$ and $x\geq
1$ both the upper and the lower bound in Theorem \ref{thm:3} as well
as the lower bound in Theorem \ref{thm:2} are asymptotically sharp
as $\lambda\rightarrow\infty$. Indeed, the extreme zeros
$x_{1}(\lambda)=-x_{n}(\lambda)$ of $P_n^{(\lambda)}$ satisfy (see,
e.g. \cite{EL:1990} or \cite{DN:2010})
\[
x_{n}^2(\lambda)=x_{1}^2(\lambda)\leq
\frac{(n-1)(n+2\lambda+1)}{(n+\lambda)^2}\rightarrow 0\ \ \mbox{ as
}\ \ \lambda\rightarrow\infty,
\]
whence
\[
\lim_{\lambda\rightarrow\infty}\frac{P_n^{(\lambda)}(x)}{a_n(\lambda)}=x^n\,,
\]
and consequently
\[
\lim_{\lambda\rightarrow\infty}u_n(x)=\frac{n}{x}\,.
\]
Clearly, the bounds for $u_n$ in Theorem \ref{thm:3} as well as the
lower bound in Theorem~\ref{thm:2} enjoy the same limit.\medskip

{\bf Acknowledgement.} The author is indebted to Professor Gancho
Tachev for kindly communicating Conjecture 1 to him. This research
was supported by the Bulgarian Science Fund through Contract no.
DDVU 02/30.

\end{document}